\documentclass[10pt,a4paper]{amsart}

\usepackage{amssymb}
\usepackage{color}
\usepackage{hyperref}
\usepackage{tikz}
\usepackage[utf8]{inputenc}

      \newtheorem{theorem}{Theorem}
[section]
      \newtheorem{lemma}[theorem]{Lemma}
      
      \newtheorem{proposition}[theorem]{Proposition}

      \newtheorem{definition}[theorem]{Definition}

\newenvironment{pf}{\begin{trivlist}\item[]{\bf Proof:\ }}
{\mbox{}\hfill\rule{.08in}{.08in}\end{trivlist}}

\newcommand{\R}{\mathbb R}
\newcommand{\Z}{\mathbb Z}

\newcommand{\Q}{\mathbb Q}
      \makeatletter
      \def\@setcopyright{}
      \def\serieslogo@{}
      \makeatother

\usepackage[arrow, matrix, curve]{xy}

\begin{document}

   \author{Sungwoon Kim}
   \address{Department of Mathematics, Jeju National University, 102 Jejudaehak-ro, Jeju, 63243, Republic of Korea}
   \email{sungwoon@jejunu.ac.kr}
   
    \author{Thilo Kuessner}
   \address{Center for Mathematical Challenges, Korea Institute for Advanced Study, Hoegiro 85, Dongdaemun-gu,
   Seoul, 130-722, Republic of Korea}
   \email{kuessner@kias.re.kr}

\title{On boundary maps of Anosov representations of a surface group to $\mathrm{SL}(3,\mathbb R)$}

\date{}

\maketitle

\begin{abstract}
We prove that Anosov representations from a closed surface group to $\mathrm{SL}(3,\R)$ are uniquely determined by their boundary maps $S^1\to \mathrm{Flag}(\R P^2)$ if and only if they do not factor over a completely reducible representation $\mathrm{SL}(2,\R)\to \mathrm{SL}(3,\R)$. 

On the other hand, there are families of completely reducible representations which can not be distinguished neither by their boundary maps nor by the topological conjugacy class of the action on their domain of discontinuity. 

We also prove that the quotient of the space of Anosov representations by the action of the mapping class group has at least $g+2$ components where $g$ is the genus of the surface.
\end{abstract}

\section{Introduction}
Let $\Sigma_g$ be an oriented, closed surface of genus $g \ge 2$. Then its fundamental group is presented as 
$$\pi_1\Sigma_g= \left\langle a_1,b_1,\ldots,a_g,b_g \left| \ \prod_{i=1}^g\left[a_i, b_i\right]=1 \right\rangle \right.$$ Here $[a, b]$ denotes the commutator of $a$ and $b$.
Given a Lie group $G$, the representation variety $\mathrm{Hom}(\pi_1\Sigma_g,G)$ has a natural topology as a subset of $G^{2g}$. When $G=\mathrm{SL}(3,\R)$, Hitchin (\cite{hit}) proved that the representation variety
$\mathrm{Hom}(\pi_1\Sigma_g,\mathrm{SL}(3,\R))$
has three connected components, two of which correspond to topologically trivial representations with vanishing Stiefel-Whitney classes. 

One of these two components is the so-called Hitchin component whose elements, the Hitchin representations, can be characterised by various equivalent properties. They are hyperconvex representations in the sense of \cite{lab} and correspond to convex projective 
structures on $\Sigma_g$ by \cite{cg}. In particular they have H\"older-embedded circles $\Lambda$ in the flag variety $\mathrm{Flag}(\R P^2)$ as limit curves. They are also characterized as 
representations with positive $X$-coordinates in the sense of \cite{fg}. This implies that all $\gamma\in\pi_1\Sigma_g$ are mapped to matrices with three distinct, positive real eigenvalues.

One important theme is that Hitchin representations are determined by and can be studied via their boundary map $\xi\colon S^1_\infty \to \mathrm{Flag}(\R P^2)$. In particular, the parametrisation of Hitchin representations by Bonahon-Dreyer (\cite{bd}) uses Fock-Goncharov's $X$-coordinates which are determined by the boundary map alone. 

We are interested in the other component of representations with trivial Stiefel-Whitney class. This component contains the trivial representation and also all representations arising from the composition of the (lift of the) monodromy of a hyperbolic structure $$\pi_1\Sigma_g\to \mathrm{SL}(2,\R)$$ with the natural reducible representation by block matrices $$\mathrm{SL}(2,\R)\to \mathrm{SL}(3,\R).$$ 
The latter examples and its deformations have been studied by Barbot \cite{bar}, who proved that for so-called radial deformations, i.e., deformations arising by multiplication with $\left(\begin{array}{ccc}e^{u(\gamma)}&0&0\\
0&e^{-2u(\gamma)}&0\\
0&0&e^{u(\gamma)}\end{array}\right)$
for some homomorphism $u\colon\pi_1\Sigma_g\to\R$, one always has the same boundary map $\xi\colon\partial_\infty \pi_1(\Sigma_g) \to \mathrm{Flag}(\R P^2)$, thus the same limit curve $\Lambda_0$
and a domain of discontinuity in the flag variety with quotient a $2$-fold covering space of the unit tangent bundle $T^1\Sigma_g$. 

In \cite{bar} Barbot asked whether the space of non-hyperconvex (i.e., non-Hitchin) representations is connected. 
Regarding this question, Thierry Barbot and 
Jaejeong Lee, in Daejeon 2014, observed that there are at least $2^{2g}+1$ connected components in the space of Anosov representations from a genus $g$ surface group to $\mathrm{SL}(3,\R)$, which gives a counterexample to the question. 

The mapping class group acts naturally on the space of Anosov representations ${\mathcal A}$ and for the action on components we have the following result.

\begin{proposition}\label{prop1.1} There are at least $g+2$ orbits for the action of the mapping class group $\mathrm{MCG}(\Sigma_g)$ on $\pi_0{\mathcal{A}}$.\end{proposition}

The disconnectedness of non-Hitchin Anosov representations is due to completely reducible Anosov representations coming 
in tuples of $2^{2g}$ representations which all have to be in distinct components of ${\mathcal A}$. These $2^{2g}$ Anosov representations all share the same boundary map $\xi\colon S^1\to \mathrm{Flag}(\R P^2)$. In particular they can 
not be distinguished in terms of invariants defined via boundary maps and in fact their Fock-Goncharov invariants are not well-defined.\footnotemark\footnotetext{So this is different from the case of representations to $\mathrm{PSL}(2,\R)$, where it is a consequence of Goldman's theorem that Anosov representations are uniquely determined by their boundary maps.} Moreover we see in \hyperref[discont]{Section \ref*{discont}} that these representations can not be distinguished by the topological actions on their domains of discontinuity.

We show however in \hyperref[unique]{Section \ref*{unique}} that this is an exceptional behaviour, i.e., that these examples (and their radial deformations as considered by Barbot) are the only Anosov representations which are not determined by their boundary maps.

\begin{theorem} If $\rho_1,\rho_2\colon\pi_1\Sigma_g\to \mathrm{SL}(3,\R)$ are Anosov representations with
the same boundary map $\xi\colon \partial_\infty \pi_1(\Sigma_g) \to \mathrm{Flag}(\R P^2)$, both factor over some completely reducible representation $\mathrm{SL}(2,\R)\to \mathrm{SL}(3,\R)$ 
and 
$\rho_2$ is obtained from $\rho_1$ by left multiplication 
with some conjugate of  
$$\left(\begin{array}{ccc}\lambda(\gamma)&0&0\\
0&\frac{1}{\lambda(\gamma)^2}&0\\
0&0&\lambda(\gamma)\end{array}\right)$$
for some homomorphism $\lambda\colon\pi_1\Sigma_g\to\R\setminus\left\{0\right\}$.\end{theorem}

In \cite[Definition 1.9]{fg}, Fock and Goncharov defined a universal higher Teichm\"uller space which in the case of $G=\mathrm{PGL}(3,\R)$ consists of 
all positive maps $\xi\colon \Q P^1\to \mathrm{Flag}(\R P^2)$ modulo the action 
of $\mathrm{PGL}(3,\R)$, and 
they showed that a subset of it parametrises the Hitchin component. Our result shows that one can still parametrise the not completely reducible Anosov representations
by a subset of the (not necessarily positive) 
maps $\xi\colon \Q P^1\to \mathrm{Flag}(\R P^2)$ modulo the action of $\mathrm{PGL}(3,\R)$ 
\\

{\bf Acknowledgements.} We thank Jaejeong Lee for raising the question answered in Proposition \ref{prop1.1}, and
Thierry Barbot for some extremely helpful discussions during the conference "Geometry of groups and spaces" 2014 in Daejeon. 

The first author gratefully acknowledges the partial support of Basic Science Research Program through the National Research Foundation of Korea funded by the Ministry of Education, Science and Technology (NRF-2015R1D1A1A09058742).

\section{Recollections}

Throughout the paper $\Sigma_g$ will be the closed, orientable surface of genus $g\ge 2$. 
We will freely use the identifications $$H^1(\Sigma_g;\Z/2\Z)\cong \mathrm{Hom}(\pi_1\Sigma_g,\Z/2\Z)\cong \mathrm{Hom}(H_1(\Sigma_g;\Z),\Z/2\Z)\cong(\Z/2\Z)^{2g}.$$
We will always denote $G=\mathrm{PGL}(3,\R)=\mathrm{SL}(3,\R)$ and $B\subset G$ will be the subgroup of upper triangular matrices.

\subsection{Anosov representations}
\subsubsection{Definitions}

Recall that a flag in $\R P^2$ is a pair 
$$(\left[v\right],\left[f\right])\in P(\R^3)\times P(\R^{3*})$$ with $f(v)=0$. Denote by $$X:=\mathrm{Flag}(\R P^2)$$ 
the flag variety of $\R P^2$ and by $$Y:=\mathrm{Frame}(\R P^2)$$ the frame variety, that is, the 
space of noncollinear triples of points in $\R P^2$. 

There 
is a well-known open embedding $\iota\colon Y\to X\times X$ 
sending $(\left[u\right],\left[v\right],\left[w\right])$ to $$((\left[u\right],\left[(uv)^*\right]),(\left[w\right],\left[(wv)^*\right]),$$ see \cite[Section 2.3]{bar}. The image of this embedding is an open set and in particular $T_yY$ is naturally identified with $T_{\iota(y)}(X\times X)$ for each $y\in Y$. The two direct summands of $T(X\times X)=TX\oplus TX$ are denoted by $E^+$ and $E^-$.

Two flags $(\left[v_1\right],\left[f_1\right])$ and $(\left[v_2\right],\left[f_2\right])$ are called transverse if $f_2(v_1)\not=0$ and $f_1(v_2)\not=0$.

Anosov representations were originally considered by Labourie in \cite{lab}, the notion 
of $P$-Anosov representations for general parabolic subgroups $P\subset G$ 
was defined by Guichard-Wienhard in \cite{gw}. In this paper we will 
only consider the case that $P=B$ is the group of upper triangular matrices and henceforth abbreviate ``$B$-Anosov representation" by ``Anosov representation". Before giving a definition of Anosov representation, a closed surface $\Sigma_g$ is assumed to be a hyperbolic surface.

\begin{definition}\label{anosov} A representation $\rho\colon\pi_1\Sigma_g\to \mathrm{SL}(3,\R)$ is an
Anosov representation if there exist continuous, $\rho$-equivariant maps $$\xi^\pm\colon\partial_\infty \mathbb H^2\to \mathrm{Flag}(\R P^2)$$ such that

\begin{itemize}
\item[(i)] $\xi^+(\eta)$ and $\xi^-(\eta)$ are transverse for each $\eta\in\partial_\infty \mathbb H^2$, so $\xi^+$ and $\xi^-$ combine to a map $$\widetilde{\sigma}\colon T^1 \mathbb H^2\to \mathrm{Frame}(\R P^2).$$

\item[(ii)] The lifted geodesic flow on $\widetilde{\sigma}^*E^+$ resp. $\widetilde{\sigma}^*E^-$ is dilating resp.\ contracting.
\end{itemize} \end{definition} 


In our case we can assume that $\xi^-=\partial_\infty s\circ\xi^+$, with $\partial_\infty s$ induced by the antidiagonal matrix $s\in \mathrm{SL}(3,\R)$ permuting the basis vectors $e_1$ and $e_3$. Thus we will
talk about ``the'' boundary map $\xi:=\xi^+$, see \cite[Section 4.5]{gw}.

\subsubsection{Space of Anosov representations}
Let us denote
$${\mathcal{A}}\subset \mathrm{Hom}(\pi_1\Sigma_g,\mathrm{SL}(3,\R))$$
the set of Anosov representations. By \cite[Proposition 2.1]{lab} it is an open subset of the representation variety. The work of Labourie shows that $\mathcal{A}$ contains the Hitchin component and the work of Barbot exhibits Anosov representations in the other component of the topologically trivial representations. By \cite[Corollary 6.6]{bar} only those two components can contain Anosov representations. Our aim is to distinguish components of $\mathcal{A}$ inside the non-hyperconvex component of topologically trivial representations in $\mathrm{Hom}(\pi_1\Sigma_g,\mathrm{SL}(3,\R))$.

Representations in this component have been studied in \cite{bar}. One of the results was that in all cases $\mathrm{Flag}(\R P^2)$ decomposes into the limit curve $\Lambda=\xi(\partial_\infty \mathbb H^2)$, a domain of discontinuity $\Omega$ homeomorphic to a solid torus, and two invariant M\"obius bands with complicated dynamics. Moreover, the quotient of $\Omega$ by the $\rho(\pi_1(\Sigma_g))$-action is a circle bundle over $\Sigma_g$.

In \cite[Section 8]{bar}, Barbot asked the following questions for non-hyperconvex Anosov representations $\rho\colon\pi_1\Sigma_g\to \mathrm{SL}(3,\R)$, noting that a positive answer to Question 2 would imply a positive answer to Question 1 (cf.\ \cite[Theorem 9.12]{gw}).\\
\\
{\bf Question 1}: Is the circle bundle $\rho(\Gamma)\backslash\Omega$ homeomorphic to the double
covering of the unit tangent bundle of $\Sigma_g$?\\
\\
{\bf Question 2}: Is the space of non-hyperconvex Anosov representations connected?\\

As mentioned before, due to the observation of T. Barbot and J. Lee (see Section \ref{redperm}), it turns out that the answer for Question 2 is no. On the other had, Question 1 might still have a positive answer in view of the result of \hyperref[discont]{Section \ref*{discont}} below.

\section{Completely reducible representations}\label{reducible}

In this section we consider the completely reducible representations $\rho_\phi$ which 
yield $2^{2g}$ different components of non-hyperconvex Anosov representations. The remainder of the section will not play a r$\hat{o}$le for this paper, though it might be of independent interest: we show that the $2^{2g}$-tuples of completely reducible representations with the same boundary map can also not be distinguished by the action on their domain of discontinuity, and we show that they are singular points of the character variety.

\subsection{Construction of a $(\Z/2\Z)^{2g}$-action}\label{red1} 

Assume a fixed representation $$\rho_0\colon\pi_1\Sigma_g\to \mathrm{SL}(2,\R)\to \mathrm{SL}(3,\R),$$ where we will as in \cite{bar} think of $\mathrm{SL}(2,\R)$ embedded in $\mathrm{SL}(3,\R)$ compatible with the embedding $(x,y)\to(x,0,y)$ of $\R^2$ in $\R^3$. Let us denote $$J_{13}=\left(\begin{array}{ccc}-1&0&0\\
0&1&0\\
0&0&-1\end{array}\right).$$
For each homomorphism $$\phi\colon\pi_1\Sigma_g\to \Z/2\Z=\left\{0,1\right\}$$ we can consider the representation
$\rho_\phi \colon\pi_1\Sigma_g\to \mathrm{SL}(3,\R)$ defined by
$$\rho_\phi(\gamma)=\rho_0(\gamma)J_{13}^{\phi(\gamma)}$$
for all $\gamma\in\pi_1\Sigma_g$. This representation is well-defined because $J_{13}$ commutes with all $\rho_0(\gamma)$ and hence the relation $\prod_{i=1}^g\left[\rho_0(a_i),\rho_0(b_i)\right]=1$ for the standard generators $a_1,b_1,\ldots,a_g,b_g$ of $\pi_1\Sigma_g$ implies $\prod_{i=1}^g\left[\rho_\phi(a_i),\rho_\phi(b_i)\right]=1$. 

Observe that even though the images in $\mathrm{SL}(2,\R)$ project to the same representations in $\mathrm{PSL}(2,\R)$, this is not the case for the images in $\mathrm{SL}(3,\R)$ in view of the equality $\mathrm{SL}(3,\R)=\mathrm{PGL}(3,\R)$.

\subsection{Domains of discontinuity}\label{discont} It is easy to check that all the $\rho_\phi$ are Anosov representations with the same boundary map as $\rho_0$, namely the embedding $\R P^1\to \mathrm{Flag}(\R P^2)$ which is induced by the embedding $\R P^1\to \R P^2$ given by $\left[x:y\right]\to\left[x:0:y\right]$. Let $L$ be the image of the latter curve in $\R P^2$, and $L^*=\left\{\left[f\right]\colon f(e_2)=0\right\}$, the image of the boundary map 
is $L\times L^*$ and one of the three components of
its complement is $$\Omega=\left\{(\left[v\right],\left[f\right])\colon v\not\in L \text{ and }f\not\in L^*\right\}\subset \mathrm{Flag}(\R P^2),$$ which can be
interpreted as the projective tangent bundle of the disk $\R P^2\setminus L$, and is thus equivariantly homeomorphic to $\mathrm{SL}(2,\R)$. The action of $\rho_\phi(\pi_1\Sigma_g)$ on $\Omega$ is properly discontinuous as a special case of \cite[Theorem 5.1]{bar}. We will argue that the actions of $\rho_\phi(\pi_1\Sigma_g)$ for different $\phi$ do all yield the same quotient manifold $\rho_\phi(\pi_1\Sigma_g)\backslash\Omega$.\\

{\bf The base space.} A hyperbolic structure on $\Sigma_g$ is given by its monodromy representation $\rho\colon\pi_1\Sigma_g\to \mathrm{PSL}(2,\R)$. The quotient $$\rho(\pi_1\Sigma_g)\backslash \mathrm{PSL}(2,\R)$$ is the unit tangent bundle $T^1\Sigma_g$. Because this is a circle bundle we have an exact sequence
$$0\to\Z\to \pi_1(T^1\Sigma_g)\to \pi_1\Sigma_g\to 1.$$
By Culler's theorem $\rho$ can be lifted to a representation $\rho_0\colon\pi_1\Sigma_g\to \mathrm{SL}(2,\R)$. We assume such a lift to be fixed.\\

{\bf The covering spaces.} Since the Euler class of $T^1\Sigma_g$ is even, the spectral sequence for the homology with $\Z/2\Z$-coefficients degenerates at the $E_2$-term (in fact the only potentially nontrivial $d_2$-differential is multiplication by the Euler class) and thus we have $$H_1(\pi_1(T^1\Sigma_g);\Z/2\Z)\simeq H_1(\Z;\Z/2\Z)\oplus H_1(\pi_1\Sigma_g;\Z/2\Z).$$ As any homomorphism from $\pi_1(T\Sigma_g)$ to $\Z/2\Z$ has to factor through the abelianization $H_1(T^1\Sigma_g)$, this implies
$$\mathrm{Hom}(\pi_1(T^1\Sigma_g),\Z/2\Z)\simeq
H^1(\Z;\Z/2\Z)\oplus H^1(\pi_1\Sigma_g;\Z/2\Z)\simeq(\Z/2\Z)^{2g+1}$$
In particular, for each homomorphism $\phi\colon\pi_1\Sigma_g\to \Z/2\Z$ we have a uniquely defined homomorphism  $$\Phi\colon\pi_1(T^1\Sigma_g)\to {\Z/2\Z}=\left\{0, 1\right\}$$ which sends the generator of $H_1(S^1;\Z/2\Z)$ to the nontrivial element $1\in\Z/2\Z$ and agrees with $\phi$ on $H_1(\Sigma_g;\Z/2\Z)$.

Inspection shows that $$\rho_\phi(\pi_1\Sigma_g)\backslash \mathrm{SL}(2,\R)$$ is the 2-fold covering space of $\rho(\pi_1\Sigma_g)\backslash \mathrm{PSL}(2,\R)$ which corresponds to the homomorphism $\Phi\colon\pi_1(T^1\Sigma_g)\to {\Z/2\Z}$.\\

{\bf Euler class.} Circle bundles over $\Sigma_g$ are classified by their Euler class $$e\in H^2(\pi_1\Sigma_g;\Z)\cong\Z.$$ It is well-known that $T^1\Sigma_g$ is a circle bundle of Euler class $2-2g$. Our quotients
$$\rho_\phi(\pi_1\Sigma_g)\backslash\Omega\cong \rho_\phi(\pi_1\Sigma_g)\backslash \mathrm{SL}(2,\R)$$
are fibre-wise double covers of $T^1\Sigma_g$ and therefore are circle bundles of Euler class $1-g$. So they are all isomorphic as circle bundles and in particular their total spaces are all homeomorphic to each other.

The homeomorphisms lift to equivariant homeomorphisms of the domains of discontinuity. So we see that the actions of the different $\rho_\phi(\pi_1\Sigma_g)$ on $\Omega$ are all topologically conjugate to each other.

\subsection{Deformations}
Although this will not be used in the remainder of the paper, we consider it worthwhile mentioning that the completely reducible representations are singular points of the character variety. Namely, it is well-known that the character variety $\mathrm{Hom}(\pi_1\Sigma_g,\mathrm{SL}(3,\R))//\mathrm{SL}(3,\R)$ has dimension $16g-16$ and we will show that at completely reducible representations the dimension of the Zariski tangent space will be $16g-14$.

We use that the dimension of the Zariski tangent space at semisimple representations is $H^1(\Gamma,Ad\circ\rho)$, see \cite{lm}. To compute the
latter we decompose the Lie algebra $\mathrm{SL}(3,\R)$ as
$$\mathrm{SL}(3,\R)=\mathrm{SL}(2,\R)\oplus\R^2\oplus\R^2\oplus\R,$$
where one $\R^2$-summand is spanned by the elementary matrices $e_{12},e_{23}$, the other $\R^2$-summand is spanned by $e_{21},e_{32}$ and the $\R$-summand is spanned by the diagonal matrix $diag(-1,2,-1)$.
The summands of this decomposition are orthogonal with respect to the Killing form and are preserved under the adjoint action $Ad$.
We note that the action of $Ad$ on the $\R$-summand is trivial. An explicit computation shows that the action of $Ad$ on the $\R^2$-summands comes from the standard linear action $\rho_{st}$ of $\mathrm{SL}(2,\R)$ on $\R^2$.
In particular, the group cohomology decomposes as a direct sum
\begin{eqnarray*} \lefteqn{H^1(\pi_1\Sigma_g,Ad\circ\iota\circ\rho_0)} \\ &=& H^1(\pi_1\Sigma_g,Ad\circ\rho_0)\oplus H^1(\pi_1\Sigma_g,\R)\oplus H^1(\pi_1\Sigma_g,\rho_{st}\circ\rho_0)\oplus H^1(\pi_1\Sigma_g,\rho_{st}\circ\rho_0) \\
 &=& T_{\rho_0}T(\Sigma_g)\oplus \R^{2g}\oplus  H^1(\pi_1\Sigma_g,\rho_{st}\circ\rho_0)\oplus  H^1(\pi_1\Sigma_g,\rho_{st}\circ\rho_0),\end{eqnarray*}
where $T(\Sigma_g)$ means Teichm\"uller space and we use that $\pi_1\Sigma_g$ acts trivially on $\R$. 

The Teichm\"uller space of $\Sigma_g$ has dimension $6g-6$ and $H^1(\Sigma_g,\R)$ has dimension $2g$. 
A variant of the Hopf trace formula argument shows 
\begin{align*} \dim\ H^0(\pi_1\Sigma_g,\R^2) & -dim\ H^1(\pi_1\Sigma_g,\R^2)+dim\ H^2(\pi_1\Sigma_g,\R^2) \\ & = \chi(\pi_1\Sigma_g)dim(\R^2)=4-4g.\end{align*} The action of the cocompact  lattice $\pi_1\Sigma_g\subset \mathrm{SL}(2,\R)$ on $\R^2$ has no nonzero fixed vector, hence $H^0(\pi_1\Sigma_g,\R^2)=(\R^2)^{\pi_1\Sigma_g}=0$ and by Poincar\'e duality $H^2(\pi_1\Sigma_g,\R^2)=0$, thus $$dim\ H^1(\pi_1\Sigma_g,\R^2)=4g-4.$$ Altogether    
$$dim\ H^1(\pi_1\Sigma_g,Ad\circ\iota\circ\rho_0)=16 g-14.$$

\section{Disconnectedness: counting components modulo the mapping class group action}\label{disconnect}

\subsection{Disconnectedness}\label{redperm}
Barbot and Lee observed that there are at least $2^{2g}$ components in the space of non-hyperconvex Anosov representations. Their proof proceeds by showing that given a completely irreducible representation $\rho_0$ as in \hyperref[red1]{Section \ref*{red1}} 
all $\rho_\phi$ with $\phi\in H^1(\Sigma_g;\Z/2\Z)$ belong to distinct path components of ${\mathcal A}$. 
For reader's convenience, we here sketch their proof.

The basic reason is that by \cite[Proposition 3.2]{lab} Anosovness of $\rho$ implies that for all $\gamma\in\pi_1\Sigma_g$ the matrix
$$\rho_\phi(\gamma)\in \mathrm{SL}(3,\R)$$
has three distinct real eigenvalues. 

So, decomposing the set of $3\times 3$-matrices with real eigenvalues into the two sets 
$$A_0=\left\{A\in \mathrm{SL}(3,\R)\colon\mbox{$A$\ has\ 3\ positive eigenvalues}\right\},$$ 
$$A_1=\left\{A\in \mathrm{SL}(3,\R)\colon\mbox{$A$\ has\ 1\ positive and 2 negative eigenvalues}\right\},$$ 
then each $\rho(\gamma)$ belongs either to $A_0$ or $A_1$ and there is an assignment
$$F\colon {\mathcal{A}}\to \mathrm{Map}(\left\{a_1,b_1,\ldots,a_g,b_g\right\},\Z/2\Z)\simeq(\Z/2\Z)^{2g}$$
by assigning for each Anosov representation $\rho\in{\mathcal A}$, each $k\in\left\{0,1\right\}$ and each of the standard generators $a_1,b_1,\ldots,a_g,b_g$ of $\pi_1\Sigma_g$
$$F(\rho)(\gamma)=k\Longleftrightarrow \rho(\gamma)\in A_k.$$
(It seems unlikely that $F(\rho)$ is a homomorphism for arbitrary $\rho\in{\mathcal A}$, although this is true for the representations $\rho_\phi$ from \hyperref[red1]{Section \ref*{red1}}.)

$F$ is surjective because each $\phi\in \mathrm{Map}(\left\{a_1,b_1,\ldots,a_g,b_g\right\},\Z/2\Z)$ is realised by 
the representation $\rho_\phi$ from \hyperref[red1]{Section \ref*{red1}}.  
On the other hand, $F$ is constant on components of ${\mathcal{A}}$ because for a continuous path $\rho_t$ of representations to $\mathrm{SL}(3,\R)$, the value of $\rho_t(\gamma)$ for some fixed $\gamma\in\pi_1\Sigma_g$ can not switch from $A_0$ to $A_1$ while $t$ is changing continuously.  This proves that ${\mathcal{A}}$ has at least $2^{2g}$ components besides the Hitchin component.

\subsection{Action of the mapping class group}\label{mcg}
The mapping class group $\mathrm{MCG}(\Sigma_g)$ (i.e., the group of diffeomorphisms modulo isotopy) of $\Sigma_g$ acts canonically on $\mathrm{Hom}(\pi_1\Sigma_g,\mathrm{SL}(3,\R))$. It obviously maps Anosov representations to Anosov representations, so we can consider 
the quotient $\mathrm{MCG}(\Sigma_g)\backslash{\mathcal{A}}$ and we are going to show that $\mathrm{MCG}(\Sigma_g)\backslash{\mathcal{A}}$ has at least $g+2$ connected components. This subsection is devoted to the proof of the following proposition.

\begin{proposition}\label{disconn2}There are at least $g+2$ orbits for the action of the mapping class group $\mathrm{MCG}(\Sigma_g)$ on $\pi_0{\mathcal{A}}$.\end{proposition}

As in \hyperref[red1]{Section \ref*{red1}} we fix a representation $\rho_0\colon \pi_1\Sigma_g\to \mathrm{SL}(2,\R)\subset \mathrm{SL}(3,\R)$ and 
will consider the finite subset ${\mathcal A}_0\subset {\mathcal A}$ which consists of the representations $\rho_\phi$ for the $2^{2g}$ different
choices of $\phi\in \mathrm{Hom}(\pi_1\Sigma_g,\Z/2\Z)$. We want to show that the elements of ${\mathcal A}_0$ belong to $g+1$ different orbits of the
maping class group. This implies the claim of \hyperref[disconn2]{Proposition \ref*{disconn2}} because the argument in \hyperref[redperm]{Section \ref*{redperm}} 
shows that all elements of ${\mathcal A}_0$ belong to pairwise distinct components of ${\mathcal A}$ (and of course the Hitchin component is preserved by the maping class group).

First we note that the action of the mapping class group on ${\mathcal{A}}_0$ and on $\mathrm{Hom}(\pi_1\Sigma_g,\Z/2\Z)$ are compatible. Indeed, for $f\in \mathrm{MCG}(\Sigma_g)$ and $\rho=\rho_\phi\in{\mathcal{A}}_0$ with 
$\phi\in \mathrm{Hom}(\pi_1\Sigma_g,\Z/2\Z)$ and for $F\in \mathrm{Hom}(\pi_1\Sigma_g,\Z/2\Z)$ we have for any $\gamma\in\pi_1\Sigma_g$
$$F(f^*\rho)(\gamma)=F(\rho)(f_*\gamma)$$
because
$F(f^*\rho)(\gamma)=0$ is equivalent to $(f^*\rho)(\gamma)$ having three positive eigenvalues, which is of course equivalent to $\rho(f_*\gamma)$ having three positive eigenvalues, hence to $ F(\rho)(f_*\gamma)=0$. 

So the orbits of the mapping class group on $\mathcal{A}_0$ are mapped to its orbits on 
$$\mathrm{Hom}(\pi_1\Sigma_g,\Z/2\Z)=H^1(\Sigma_g,\Z/2\Z).$$
As noted in \hyperref[redperm]{Section \ref*{redperm}} every $\phi\in \mathrm{Hom}(\pi_1\Sigma_g,\Z/2\Z)$ is realised by the representation $\rho_\phi
\in{\mathcal A}_0$ and thus \hyperref[disconn2]{Proposition \ref*{disconn2}} is a consequence of the following lemma.

\begin{lemma}\label{mcgaction} There are $g+1$ orbits for the action of $\mathrm{MCG}(\Sigma_g)$ on $\mathrm{Hom}(\pi_1\Sigma_g,\Z/2\Z)$.\end{lemma}
\begin{pf}
It is well-known (and easy to prove) that the {\em intersection form modulo 2} 
$$i\colon H_1(\Sigma_g,\Z/2\Z)\times H_1(\Sigma_g,\Z/2\Z)\to \Z/2\Z$$
defines a symplectic form on the ${\Z/2\Z}$-vector space $ 
H_1(\Sigma_g,\Z/2\Z)$ and that the action of the mapping class group preserves this symplectic form. 

For each $$\phi\in \mathrm{Hom}(\pi_1\Sigma_g,\Z/2\Z)=\mathrm{Hom}(H_1(\Sigma_g,\Z/2\Z),\Z/2\Z)$$
we let $d_\phi$ be the dimension of the maximal symplectic subspace on which $\phi$ is constant $0$. This number is invariant under the action of the mapping class group, thus $$d_{\phi_1}\not=d_{\phi_2}$$ implies that $\phi_1$ and $\phi_2$ are not in the same $\mathrm{MCG}(\Sigma_g)$-orbit. 

The number $d_\phi$ can take integer even values in 
$$\left\{0,2,4,\ldots,2g\right\}.$$ 
We claim that each of these $g+1$ values can indeed 
be realised for some $\phi$.
An explicit realisation for a given $d$ is for example given as follows. Let $$a_1,b_1,\ldots,a_g,b_g$$ be the standard basis of $H_1(\Sigma_g,\Z/2\Z)$ with respect to which the intersection form is given by the standard symplectic form $$i(a_k,b_l)=\delta_{kl},i(a_k,a_l)=i(b_k,b_l)=0$$ for $k,l=1,\ldots,g$. Then, for a given $d\in\left\{0,2,4,\ldots,2g\right\}$, the homomorphism $$\phi_d(a_1)=\phi_d(b_1)=\ldots=\phi(a_d)=\phi(b_d)=0,$$
$$ \phi_d(a_{d+1})=\phi_d(b_{d+1})=\ldots=\phi_d(a_g)=\phi_d(b_g)=1$$ obviously realises $d_\phi=d$.

This part of \hyperref[mcgaction]{Lemma \ref*{mcgaction}} actually suffices to prove \hyperref[disconn2]{Proposition \ref*{disconn2}} but for 
completeness we still show that there are {\em exactly} $g+1$ orbits for the action of the mapping class group on 
$\mathrm{Hom}(\pi_1\Sigma_g,\Z/2\Z)$ (and thus that our completely reducible examples yield {\em exactly} g+1 distinct components of $\mathrm{MCG}(\Sigma_g)\backslash{\mathcal{A}}$). 
For this we have to show that any $\phi\in \mathrm{Hom}(\pi_1\Sigma_g,\Z/2\Z)$ is in the $\mathrm{MCG}(\Sigma_g)$-orbit of $\phi_{d_\phi}$, where $d_\phi$ is 
the dimension of the maximal symplectic subspace on which $\phi$ is constant $0$ and $\phi_{d_\phi}$ is defined as in the previous 
paragraph. 

By its definition, the symplectic structure $i$ is standard, i.e. decomposes into the 2-dimensional subspaces generated by $a_k,b_k$ for $k=1,\ldots,g$. 
If $$\phi\not=\phi_{d_\phi},$$ then there is some set of 2-dimensional subspaces on which $\phi$ does not agree with $\phi_d$ but is not constant $0$. Say these subspaces are generated by $$a_{k_1},b_{k_1},\ldots,a_{k_l},b_{k_l}.$$ 
One can find some mapping class which sends $a_{k_1}$ to $a_1$, $b_{k_1}$ to $b_1$, ..., $a_{k_l}$ to $a_l$, $b_{k_l}$ to $b_l$, so we can w.l.o.g. assume $k_1=1,\ldots,k_l=l$. (The existence of such a mapping class follows either from the Dehn-Nielsen Theorem \cite{nie}, which says that any automorphism of a surface group is induced by some mapping class, or in this special case also from an explicit construction.)

For each $k\in\left\{1,\ldots,l\right\}$ we can apply a Dehn twist 
to map the restriction of $\phi$ to the restriction of $\phi_{d_\phi}$. In formula: assume w.l.o.g. $$\phi(a_k)=0,\phi(b_k)=1.$$ The Dehn twist $D_k:=D_{b_k}$ 
at $b_k$ sends $a_k$ to $a_k+b_k$ and $b_k$ to $b_k$. Thus 
$$\phi(D_{b_k}(a_k))=\phi(D_{b_k}(b_k))=1,$$ 
which means $$\phi(D_{b_k}(.))=\phi_{d_\phi}(.)$$ on 
the subspace generated by $a_k,b_k$.
Composition of the Dehn twists $D_{1},\ldots,D_{l}$ then yields the wanted result:
$$\phi(D_{l}\ldots D_{1}(h))=\phi_{d_\phi}(h)$$
for all $h\in H_1(\Sigma_g,\Z/2\Z).$
\end{pf}

Remark: The same argument can be used to show that the lift of the Teichm\"uller space to $\mathrm{Hom}(\pi_1\Sigma_g,\mathrm{SL}(3,\R))$ falls into $g+1$ components modulo the action of the mapping class groups.

\section{Boundary maps determine Anosov representations - almost always}\label{unique}
In \hyperref[reducible]{Section \ref*{reducible}} we have seen distinct 
Anosov representations with the same boundary map. In this section we will see that these examples are essentially the only ones for which the boundary map does not determine the Anosov representation. 
\begin{theorem}\label{uniqueness} Let $$\rho_1,\rho_2\colon\pi_1\Sigma_g\to \mathrm{SL}(3,\R)$$ be Anosov representations such that there exists a map $$\xi\colon \partial_\infty \mathbb H^2\to \mathrm{Flag}(\R P^2)$$ satisfying the conditions of (the remark after) \hyperref[anosov]{Definition \ref*{anosov}} and equivariant for both, $\rho_1$ and $\rho_2$. Then either $\rho_1=\rho_2$ or there exists a representation $r\colon \pi_1\Sigma_g\to \mathrm{SL}(2,\R)$, a homomorphism $\lambda\colon\pi_1\Sigma_g\to 
\R\setminus\left\{0\right\}$ and some $A\in \mathrm{PGL}(3,\R)$ such that
$$A\rho_1A^{-1}=\iota\circ r$$
$$A\rho_2A^{-1}=\left(\begin{array}{ccc}\lambda&0&0\\
0&\frac{1}{\lambda^2}&0\\
0&0&\lambda\end{array}\right) \circ \iota \circ r,$$
where $\iota\colon \mathrm{SL}(2,\R)\to \mathrm{SL}(3,\R)=\mathrm{PGL}(3,\R)$ is the completely reducible representation from \hyperref[reducible]{Section \ref*{reducible}}.
\end{theorem}

\begin{pf}
Consider an ideal triangulation $\Upsilon$ of $\Sigma_g$, let $\widetilde{\Upsilon}$ be the lifted ideal triangulation of $\mathbb H^2$ and (after identifying $\widetilde{\Sigma_g}$ to $\mathbb H^2$ via some hyperbolic metric) $\widetilde{\Upsilon}_0\subset\partial_\infty \mathbb H^2$ its $0$-skeleton. 

Assume $\rho_1\not=\rho_2$, so there is some $\gamma_0\in\pi_1\Sigma_g$ with $\rho_1(\gamma_0)\not=\rho_2(\gamma_0)$. For every vertex $v\in \widetilde{\Upsilon}_0$ we have $\rho_1(\gamma_0)\xi(v)=\xi(\gamma_0 v)=\rho_2(\gamma_0)\xi(v)$ and thus $$\rho_1^{-1}(\gamma_0)\rho_2(\gamma_0)\in \mathrm{Stab}(\xi(v)).$$
In particular, for every ideal triangle $T=(v_0,v_1,v_2)\in\widetilde{\Upsilon}$, $\rho_1^{-1}(\gamma_0)\rho_2(\gamma_0)$ stabilizes the associated triple of flags $(\xi(v_0),\xi(v_1),\xi(v_2))$.

By part i) of \hyperref[anosov]{Definition \ref*{anosov}} we have that $\xi(v)$ and $\xi(w)$ are transverse for all $v,w\in \widetilde{\Upsilon}_0$. An elementary argument, given in \cite[Section 2]{kue} shows that every triple of pairwise transverse flags in $\R P^2$ is in the $\mathrm{PGL}(3,\R)$-orbit of one of the following triples:
$$\left\{\left((e_1,e_3^\perp),(e_3,e_1^\perp),(e_1-e_2+e_3,(e_1+(1+X)e_2+Xe_3)^\perp\right), X\in\R\setminus\left\{0,-1\right\}\right\}$$
$$((e_1,e_3^\perp),(e_3,e_1^\perp),(e_1+e_3,(e_1+e_2-e_3)^\perp))$$
$$((e_1,e_3^\perp),(e_3,e_1^\perp),(e_1+e_2+e_3,(e_1-e_3)^\perp))$$
$$(e_1,e_3^\perp),(e_3,e_1^\perp),(e_1+e_3,(e_1-e_3)^\perp)$$

One easily checks that the last triple is the only one of these possibilities which has a nontrivial stabilizer in $\mathrm{PGL}(3,\R)$. In fact, the stabilizer of a triple in the $\mathrm{PGL}(3,\R)$-orbit of $(e_1,e_3^\perp),(e_3,e_1^\perp),(e_1+e_3,e_1^\perp-e_3^\perp))$ is conjugate to $$\left\{\left(\begin{array}{ccc}\lambda&0&0\\
0&\frac{1}{\lambda^2}&0\\
0&0&\lambda\end{array}\right) : \lambda\in\R\setminus\left\{0\right\}\right\}.$$

So $\rho_1^{-1}(\gamma_0)\rho_2(\gamma_0)\in \mathrm{Stab}(\xi(v))\setminus\left\{id\right\}$ for every $v\in \widetilde{\Upsilon}_0$ implies that for every ideal triangle $T=(v_0,v_1,v_2)\in\widetilde{\Upsilon}$ the associated triple of flags $$(\xi(v_0),\xi(v_1),\xi(v_2))$$ must be in the $\mathrm{PGL}(3,\R)$-orbit of $$((e_1,e_3^\perp),(e_3,e_1^\perp),(e_1+e_3,(e_1-e_3)^\perp).$$ 
In other words, if we denote $\xi(v_i)=(p_i,l_i)$ for $i=0,1,2$, then $p_2$ must be on the line through $p_0$ and $p_1$ and $l_2$ goes through the intersection point of $l_0$ and $l_1$.  

We claim that this implies that all $\xi(v), v\in\widetilde{\Upsilon}_0$ are of the form 
$\xi(v)=(p_v,l_v)$ with all $p_v$ lying on the same line $l$, and all $l_v$ intersecting in the same point $p$. To see this. let $\lambda^{closed}\subset \Sigma_g$ be 
the union of closed leaves of $\Upsilon$, let $\tilde{\lambda}^{closed}\subset \mathbb H^2$ its lift to $\widetilde{\Sigma}_g\simeq 
\mathbb H^2$ and let $\widetilde{U}$ be the 
component of $\mathbb H^2\setminus\tilde{\lambda}^{closed}$ containing $T$. 
Looking at the dual tree of the ideal triangulation $\widetilde{\Upsilon}\mid_{\widetilde{U}}$ we can enumerate 
its triangles such that each $T_k$ is adjacent to exactly one triangle $T_j$ 
with $j<k$, see the proof of \cite[Lemma 21]{bd}. Using this enumeration we see by induction that the claim is 
true for all vertices of triangles in $\widetilde{U}$. Further this 
also applies to the leaves of $\tilde{\lambda}^{closed}$ adjacent to $\widetilde{U}$, because one of their ideal vertices is actually a vertex of $\widetilde{\Upsilon}\mid_{\widetilde{U}}$, while the other ideal vertex is an accumulation point of 
vertices and so the claim holds by continuity of $\xi$. Next we can extend 
this argument to the components of $\mathbb H^2\setminus\tilde{\lambda}^{closed}$ 
adjacent to $\widetilde{U}$. Namely, the same argument shows the claim (with a priori possibly different $(p,l)$) 
for all triangles in an adjacent component. But since the $(p,l)$ agree on both ideal vertices of the leaf in $\tilde{\lambda}^{closed}$ 
along which they are adjacent, the $(p,l)$ 
must actually be the same for both adjacent components, just because there is a unique line through two points and a unique point common to two lines. 
Finally we use the dual tree to the decomposition into components of 
$\mathbb H^2\setminus\tilde{\lambda}^{closed}$ to enumerate these components 
such that each $U_k$ is adjacent to exactly 
one component $U_j$ with $j<k$, as in the proof of \cite[Lemma 24]{bd}, so that we can induct on the components and finally get the claim for all triangles in all components.

Now fix an
ideal triangle $T=(v_0,v_1,v_2)\in\widetilde{\Upsilon}$ and a projective map $A\in \mathrm{PGL}(3,\R)$ that sends
$$(\xi(v_0),\xi(v_1),\xi(v_2))=((p_0,l_0),(p_1,l_1),(p_2,l_2))$$ to 
$$((e_1,e_3^\perp),(e_3,e_1^\perp),(e_1+e_3,e_1^\perp-e_3^\perp)).$$ The map $$A\xi\colon\partial_\infty \mathbb H^2\to \mathrm{Flag}(\R P^2)$$ 
is equivariant for $A\rho_1A^{-1}$ and $A\rho_2A^{-1}$, and still satisfies the conditions from \hyperref[anosov]{Definition \ref*{anosov}}.

We note that $$Ap=Al_0\cap Al_1\cap Al_2=e_3^\perp\cap e_1^\perp\cap (e_1-e_3)^\perp=e_2$$
and that $Al$ is the line containing $Ap_0=e_1, Ap_1=e_3$ and $Ap_2=e_1-e_3$ and thus 
$$Al=e_2^\perp.$$

For each $\gamma\in \pi_1\Sigma_g$, we have (because $\widetilde{\Upsilon}$ is defined by lifting $\Upsilon$) that $\gamma T=(\gamma v_0,\gamma v_1,\gamma v_2)$ is one of the ideal 
triangles of $\widetilde{\Upsilon}$. 
Then $A\rho_1(\gamma)A^{-1}e_1$ is the point component of $$A\rho_1(\gamma)A^{-1}(A\xi(v_0))=A\xi(\gamma  v_0)$$ and thus belongs to $e_2^\perp$. In the same way,  $A\rho_1(\gamma)A^{-1}e_3$ is the point component of $$A\rho_1(\gamma)A^{-1}(A\xi(v_1))=A\xi(\gamma  v_1)$$ and thus belongs to $e_2^\perp$. This shows that $A\rho_1(\gamma)A^{-1}$ sends $e_2^\perp$ to itself. Similarly, considering the line components, we can show that 
$A\rho_1(\gamma)A^{-1}$ sends $e_2$ to itself. Since this is true for any $\gamma\in \pi_1\Sigma_g$ we have that the image of $A\rho_1A^{-1}$ is 
in the image of the completely reducible representation $$\iota\colon \mathrm{SL}(2,\R)\to \mathrm{SL}(3,\R)=\mathrm{PGL}(3,\R).$$ The same argument applies to $A\rho_2A^{-1}$.

Let $$r\colon \pi_1\Sigma_g\to \mathrm{SL}(2,\R)$$ be the representation over which $A\rho_1A^{-1}$ factors. Then, for each $\gamma\in\pi_1\Sigma_g$, both $A\rho_1(\gamma)A^{-1}$ and $A\rho_2(\gamma)A^{-1}$ send $A\xi(T)$ to $A\xi(\gamma T)$ and thus
$A\rho_1^{-1}(\gamma)\rho_2(\gamma)A^{-1}$ is in the stabilizer of $$A\xi(T)=((e_1,e_3^\perp),(e_3,e_1^\perp),(e_1+e_3,(e_1-e_3)^\perp)).$$ So there is a unique $\lambda(\gamma)\in\R\setminus\left\{0\right\}$ with 
$$A\rho_1^{-1}(\gamma)\rho_2(\gamma)A^{-1}=\left(\begin{array}{ccc}\lambda(\gamma)&0&0\\
0&\frac{1}{\lambda(\gamma)^2}&0\\
0&0&\lambda(\gamma)\end{array}\right),$$
from which the claimed formula for $A\rho_2(\gamma)A^{-1}$ follows. 

From the fact that the diagonal matrices of the form $\mathrm{Diag}(\lambda,\frac{1}{\lambda^2},\lambda)$ commute with the image of the completely reducible representation $\iota$, one easily concludes that $\lambda$ is a homomorphism. Indeed, this follows from the computation
\begin{eqnarray*} \lefteqn{\left(\begin{array}{ccc}\lambda(\gamma_1\gamma_2)&0&0\\
0&\frac{1}{\lambda(\gamma_1\gamma_2)^2}&0\\
0&0&\lambda(\gamma_1\gamma_2)\end{array}\right)  (\iota \circ r(\gamma_1\gamma_2))} \\ 
&=& A\rho_2(\gamma_1\gamma_2)A^{-1} \\
&=& A\rho_2(\gamma_1)A^{-1}A\rho_2(\gamma_2)A^{-1} \\
&=& \left(\begin{array}{ccc}\lambda(\gamma_1)&0&0\\
0&\frac{1}{\lambda(\gamma_1)^2}&0\\
0&0&\lambda(\gamma_1)\end{array}\right)  (\iota \circ r(\gamma_1))\left(\begin{array}{ccc}\lambda(\gamma_2)&0&0\\
0&\frac{1}{\lambda(\gamma_2)^2}&0\\
0&0&\lambda(\gamma_2)\end{array}\right)  (\iota \circ r(\gamma_2)) \\
&=& \left(\begin{array}{ccc}\lambda(\gamma_1)&0&0\\
0&\frac{1}{\lambda(\gamma_1)^2}&0\\
0&0&\lambda(\gamma_1)\end{array}\right)  
\left(\begin{array}{ccc}\lambda(\gamma_2)&0&0\\
0&\frac{1}{\lambda(\gamma_2)^2}&0\\
0&0&\lambda(\gamma_2)\end{array}\right)  (\iota \circ r(\gamma_1\gamma_2)),\end{eqnarray*}
where the last equality uses the fact that $\mathrm{Diag}(\lambda(\gamma_2),\frac{1}{\lambda(\gamma_2)^2},\lambda(\gamma_2))$ commutes
with $\iota\circ r(\gamma_1)$ and that $\iota\circ r$ is a homomorphism.
\end{pf}\

It is perhaps worth mentioning that, given a completely reducible Anosov representation $\rho_1$, not every homomorphism $\lambda$ will yield an Anosov representation $\rho_2$ as in \hyperref[uniqueness] {Theorem \ref*{uniqueness}}. Actually, \cite[Theorem 4.2]{bar} gives a precise condition for the stable norm of $\log(\lambda)$ to guarantuee that $\rho_2$ also is Anosov.

Let us finally mention that the assumption of \hyperref[uniqueness] {Theorem \ref*{uniqueness}} can not be weakened to consider only the boundary map with image in $\R P^2$ or its dual, rather than in $\mathrm{Flag}(\R P^2)$. Barbot constructs in 
\cite[Lemma 4.5]{bar} a family of (reducible, but not completely reducible) representations which all have the same boundary map to the dual space of $\R P^2$, however their boundary maps to $\R P^2$ do not agree. Similarly, one can construct a family of representations with the same boundary map to $\R P^2$, but their boundary maps to the dual space and hence to 
$\mathrm{Flag}(\R P^2)$ will not agree. 

These examples are not Zariski-dense and so they show that the class of representations determined by their boundary map is 
strictly larger than the set of Zariski-dense representations.


\begin{thebibliography}{99}
\bibitem{bar}T.\ Barbot: "Three-dimensional Anosov flag manifolds." Geom.\ Topol.\ 14 (2010), 153-191. 
\bibitem{bd}F.\ Bonahon, G.\ Dreyer: "Parameterizing Hitchin components." Duke Math.\ J.\ 163 (2014), 2935-2975.
\bibitem{cg}S.\ Choi, W.\ M.\ Goldman: "Convex real projective structures on closed surfaces are closed." Proc. Amer. Math. Soc. 118 (1993), no. 2, 657-661.
\bibitem{fm}B.\ Farb, D.\ Margalit: "A primer on mapping class groups." Princeton Mathematical Series, 49. Princeton University Press, Princeton, NJ, 2012.
\bibitem{fg}V.\ Fock, A.\ Goncharov: "Moduli spaces of local systems and higher Teichm\"uller theory." Publ.\ Math.\ Inst.\ Hautes \`Etudes Sci.\ 103 (2006), 1-211.
\bibitem{gw}O.\ Guichard, A.\ Wienhard: "Anosov representations: domains of discontinuity and applications." Inv.\ Math.\ 190 (2012), 357-438.
\bibitem{hit}N.\ J.\ Hitchin: "Lie groups and Teichm\"uller space." Topology 31 (1992), 449-473.
\bibitem{kue}T.\ Kuessner: ``Non-Hausdorff parameters for Anosov representations to SL(3,R)'', in preparation.
\bibitem{lab}F.\ Labourie: "Anosov flows, surface groups and curves in projective space." Inv.\ Math.\ 165 (2006), 51-114.
\bibitem{lm}A.\ Lubotzky, A.\ Magid: "Varieties of representations of finitely generated groups." Mem.\ AMS 58 (1985), no. 336.
\bibitem{nie}J.\ Nielsen: ''Untersuchungen zur Topologie der geschlossenen zweiseitigen Flächen.'' Acta Math.\ 50 (1927), 189-358.

\end{thebibliography}
\end{document}